\documentclass{article}

\usepackage{arxiv}

\usepackage[utf8]{inputenc}
\usepackage[T1]{fontenc}
\usepackage{hyperref}
\usepackage{url}
\usepackage{booktabs}

\usepackage{amsmath,amssymb,eucal}

\usepackage{amsfonts}

\title{Strong Converse Inequalities for Bernstein Polynomials with Explicit Asymptotic  Constants}

\newif\ifuniqueAffiliation

\uniqueAffiliationtrue

\ifuniqueAffiliation
\author{Jos\'e A. Adell  \\
	Departamento de M\'etodos Estad\'\i sticos\\
	Universidad de Zaragoza\\
	50009 Zaragoza, Spain \\
	\texttt{adell@unizar.es} \\
	\And
	Daniel C\'ardenas-Morales \\
	Departamento de Matem\'aticas\\
	Universidad de Ja\'en\\
	23071 Ja\'en, Spain \\
	\texttt{cardenas@ujaen.es} \\
}
\else

\fi

\renewcommand{\shorttitle}{}

\renewcommand{\headeright}{}

\renewcommand{\undertitle}{}

\hypersetup{pdftitle={},pdfsubject={},pdfauthor={},pdfkeywords={},}

\newtheorem{theorem}{Theorem}
\newtheorem{lemma}[theorem]{Lemma}
\newtheorem{proposition}[theorem]{Proposition}

\begin{document}

\maketitle

\begin{abstract}
We obtain strong converse inequalities for the Bernstein polynomials with explicit asymptotic constants. We give different estimation procedures in the central and non-central regions of $[0,1]$. The main ingredients in our approach are the following: representation of the derivatives of the Bernstein polynomials in terms of the Krawtchouk polynomials, estimates of different inverse moments of various random variables, sharp estimates of both absolute central moments of Bernstein polynomials and the total variation distance between binomial and Poisson distributions, and iterates of the Bernstein polynomials, together with their probabilistic representations.
\end{abstract}



\section{Introduction and statement of the main results}\label{sec1}

The Bernstein polynomials represent the most paradigmatic example of positive linear operators. Recall that for a function $f:[0,1]\rightarrow \mathbb{R}$ and a natural number $n\in \mathbb{N}$, the $n$th Bernstein polynomial of $f$ is defined as
\begin{equation}\label{1}
B_nf(x)=\sum_{k=0}^n
f\left(\frac{k}{n}\right)\binom{n}{k}x^k(1-x)^{n-k},\quad x\in [0,1].
\end{equation}
In the nineties of the last century, Ditzian and Ivanov \cite{ditzianivanov} and Totik \cite{totik} characterized the rates of uniform convergence of sequences of positive linear operators $L_nf$ towards the function $f$, as $n$ tends to infinity, in terms of the so called second order Ditzian-Totik modulus of smoothness of $f$ (cf. Ditzian and Totik \cite{ditziantotik}). In the case of Bernstein polynomials, it turns out that
\begin{equation}\label{2}
K_1\omega_2^{\varphi}\left(f;\frac1{\sqrt{n}}\right)\leq \|B_nf-f\|\leq K_2\omega_2^{\varphi}\left(f;\frac1{\sqrt{n}}\right),\quad f\in C[0,1],\ n\in \mathbb{N},
\end{equation}
for some absolute positive constants $K_1$ and $K_2$. We recall here  some necessary definitions and notations. As usual, we denote by $C[0,1]$ the space of all real continuous functions defined on $[0,1]$ endowed with the supremum norm $\|\cdot \|$, and by $C^m[0,1]$ the subspace of all $m$-times continuously differentiable functions. Besides, $f^{(m)}$ and $f^m$ represent, respectively, the $m$th derivative of $f$, and $f$ raised to the power of $m$. The second order central difference of $f$ is given by
\[
\Delta_h^2f(x)=f(x+h)-2f(x)+f(x-h),\quad h\geq 0,\quad x\pm h\in [0,1],
\]
and the Ditzian-Totik modulus of smoothness of $f$ with weight function
\[
\varphi (x)=\sqrt{x(1-x)},\quad 0\leq x\leq 1,
\]
is defined as
\[
\omega_{2}^{\varphi}(f;\delta )=\sup \left\{  \left\lvert \Delta_{h\varphi (x)}^{2}f(x) \right\rvert  :\ 0\leq h\leq \delta,\ x\pm h\varphi (x)\in [0,1] \right\},\quad \delta \geq 0.
\]

The second inequality in (\ref{2}) is called a direct inequality. Different authors completed this inequality by showing specific values for the constant $K_2$. In this regard, Gavrea et al. \cite{gavreaetall} and Bustamante \cite{bustamante} provided $K_2=3$, whereas P\u{a}lt\u{a}nea \cite{paltanea2004} obtained $K_2=2.5$, which is the best result up to date and up to our knowledge (see also \cite{rm2022} for an asymptotic result). It was also observed that if additional smoothness conditions on $f$ are added, then a direct inequality  may be valid for values of $K_2$ smaller than $2.5$. In this respect, the authors recently proved \cite{rm2022} that for any $f\in C^2[0,1]$ and $n\in \mathbb{N}$, the following inequality holds true
\begin{equation}\label{3}
\left\lvert \|B_nf-f\|-\frac{1}{2}\omega_2^{\varphi}\left(f;\frac{1}{\sqrt{n}}\right)\right \rvert \leq \frac{1}{4n}\left(\omega_1\left(f^{(2)};\frac{1}{3\sqrt{n}}\right)+\frac{1}{4}\omega_2^{\varphi}\left(f^{(2)};\frac{1}{\sqrt{n}}\right)\right),
\end{equation}
where $\omega_1 (f^{(2)};\cdot )$ stands for the usual first modulus of continuity of $f^{(2)}$. This inequality completes asymptotic results previously shown by Bustamante and Quesada \cite{bustamantequesada} and P\u{a}lt\u{a}nea \cite{paltanea2018}.

The first inequality in (\ref{2}) is called strong converse inequality, according to the terminology given in Ditzian and Ivanov \cite{ditzianivanov}, and it turns out to be the best possible one among different types of converse inequalities. It is a consequence of the works by Totik \cite{totikk}, Knopp and Zhou \cite{knoopzhou1994,knoopzhou1995}, and Sang\"{u}esa \cite{sanguesa} for general sequences of positive linear operators.

Direct and converse inequalities in the $L_p$-norm, $1\leq p\leq \infty$, have been widely considered in the literature (see, for instance, Totik \cite{totik,totikk}, Ditzian and Ivanov \cite{ditzianivanov}, Chen and Ditzian \cite{chenditzian}, Della Vecchia \cite{dellavecchia}, Guo and Qi \cite{guoqi}, Finta \cite{finta}, Gadjev \cite{gadjev}, and Bustamante \cite{bustamantee}, among many others). The main tool to prove converse inequalities is the use of $K$-
functionals, shown to be equivalent to the corresponding Ditzian-Totik modulus of smoothness. In some papers (cf. Knopp and Zhou \cite{knoopzhou1994,knoopzhou1995} and Sang\"{u}esa \cite{sanguesa}), the authors use an arbitrary number of iterates of the operators under consideration as the main tool. In any case, the proofs of converse inequalities are very involved and do not provide explicit constants in general.

This is the motivation of the present paper, whose main result reads as follows.

\begin{theorem}\label{th1}
Let $\mathfrak{F}$ denote the set of all non-affine functions $f\in C[0,1]$. Then,
\[
5\leq \varlimsup _{n\rightarrow \infty}\sup_{f\in \mathfrak{F}}\frac{\omega_2^{\varphi}\left(f;1/\sqrt{n}\right)}{\|B_nf-f\|}< 75.
\]
\end{theorem}

An analogous result to Theorem \ref{th1} referring to a non-centered gamma-type operator was given in \cite{adellsanguesa}. The upper constant there was $105$, instead of $75$. On the other hand, as seen in (\ref{3}), the estimate of the constant $K_2$ is improved when we restrict our attention to functions in $C^2[0,1]$. A similar property holds concerning strong converse inequalities, as the following result shows. In this regard, denote by $\omega_2(f;\cdot )$ the usual second modulus of continuity of $f$. Let $\mathfrak{G}\subseteq \mathfrak{F}$ be the subset of those functions $f$ such that
\begin{equation}\label{AAAAA}
    \lim _{\delta \to 0}\frac{\omega_2(f;\delta)}{\delta}=0.
\end{equation}

\begin{theorem}\label{th2}
We have
\begin{equation}\label{BBBBB}
\varlimsup _{n\rightarrow \infty}\sup_{f\in \mathfrak{G}}\frac{\omega_2^{\varphi}\left(f;1/\sqrt{n}\right)}{\|B_nf-f\|}\leq 4+ \frac{\sqrt{2}\left(\sqrt{2}+1\right)}{1-0.99/\sqrt{3}}\log 4=15.0477\ldots
\end{equation}
Moreover, $\mathfrak{F}\cap C^1[0,1]\subseteq \mathfrak{G}$.
\end{theorem}

To prove the aforementioned results, we follow an approach inspired by the work of Sang\"{u}esa \cite{sanguesa}, using different estimation procedures in the central and non-central regions of $[0,1]$. In the central region, the starting point is a Taylor's formula of third order for the third iterate of the Bernstein polynomials. As specific tools in this region, we use sharp estimates for the first absolute central moments of the Bernstein polynomials, as well as for the total variation distance between binomial and Poisson distributions. In the non-central region, we start from a formula that involves the second derivatives of the iterates of the Bernstein polynomials, providing at the same time probabilistic representations of such iterates.

Two common tools applied in both regions should be emphasized. On the one hand, several expressions of the derivatives of the Bernstein polynomials, particularly that in terms of the orthogonal polynomials with respect to the binomial distribution, namely, the Krawtchouk polynomials. On the other hand, accurate estimates of various inverse moments involving different random variables.

The main auxiliary results to give estimates in the central (resp. non-central) region are Theorems \ref{th3}, \ref{th4}, and \ref{th14} (resp. Theorems \ref{th3} and \ref{th19}). On the other hand, in some stages of the proofs of Theorems \ref{th1} and \ref{th2}, we use numerical computations performed with the software Mathematica 10.2 in order to obtain explicit constants. We finally mention that the problem of giving non-asymptotic estimates for $K_1$ remains open.

\section{General auxiliary results}\label{sec2}
Let $\mathbb{N}_0=\mathbb{N}\cup \{0\}$. From now on, whenever we write $f$,  $x$, and $n$, we respectively assume that $f\in C[0,1]$, $x\in (0,1)$, and $n\in \mathbb{N}$ is large enough.

Let $\left(\beta_m\right)_{m\geq 1}$ be a sequence of random variables such that $\beta_m$ has the beta density
\begin{equation}\label{uno}
\rho _m(\theta )=m(1-\theta)^{m-1},\quad 0\leq \theta \leq 1.
\end{equation}
These random variables allow us to write the remainder of Taylor's formula in closed form. In fact, if $g\in C^m[0,1]$, $m\in \mathbb{N}$, we have
\begin{equation}\label{dos}
g(y)=\sum_{j=0}^{m-1}\frac{g^{(j)}(x)}{j!} (y-x)^j+\frac{(y-x)^m}{m!}\mathbb{E}g^{(m)}\left(x+(y-x)\beta_m\right),\quad 0\leq y \leq 1,
\end{equation}
where $\mathbb{E}$ stands for mathematical expectation.

\begin{theorem}\label{th3}
We have
\[
\omega_2^{\varphi}\left(f;\frac{1}{\sqrt{n}}\right)\leq 4\|B_nf-f\|+\frac{\log 4}{n}\left\|\varphi ^2(B_nf)^{(2)}\right\|.
\]
\end{theorem}

\textit{Proof}.\
Suppose that $[x-\varphi (x)h,x+\varphi (x)h]\subseteq [0,1]$, with $0\leq h\leq 1/\sqrt{n}$. This implies that
\begin{equation}\label{tres}
r(x):=\frac{\varphi (x)h}{x}\in [0,1],\qquad s(x):=\frac{\varphi (x)h}{1-x}\in [0,1].
\end{equation}
We claim that
\begin{equation}\label{cuatro}
\left\vert \Delta ^2_{h\varphi (x)}g(x)\right\vert \leq h^2\left\|\varphi ^2 g^{(2)}\right\|\log 4,\quad g\in C^2[0,1].
\end{equation}
Actually, we have from (\ref{dos})
\begin{align}\label{cinco}
 \left|\Delta ^2_{h\varphi (x)}g(x)\right| &=\frac{h^2 \varphi ^2(x)}{2}\left|\mathbb{E}g^{(2)}(x-\varphi (x)h\beta_2)+\mathbb{E}g^{(2)}(x+\varphi (x)h\beta_2)\right| \cr
 &\leq h^2\left\|\varphi ^2 g^{(2)}\right\|\frac{\varphi ^2(x)}{2}\left(\mathbb{E}\frac{1}{\varphi^2\left(x-\varphi(x)h\beta_2\right)}+\mathbb{E}\frac{1}{\varphi^2\left(x+\varphi(x)h\beta_2\right)}\right).
\end{align}
Since
\begin{equation}\label{seis}
\frac{1}{\varphi^2(y)}=\frac{1}{y}+\frac{1}{1-y},\quad y\in (0,1),
\end{equation}
we get from (\ref{tres})
\[
\frac{\varphi^2(x)}{2}\left(\mathbb{E}\frac{1}{\varphi^2\left(x-\varphi(x)h\beta_2\right)}+\mathbb{E}\frac{1}{\varphi^2\left(x+\varphi(x)h\beta_2\right)}\right)
\]
\begin{equation}\label{siete}
=\frac{1-x}{2}\mathbb{E}\left(\frac{1}{1-r(x)\beta_2}+\frac{1}{1+r(x)\beta_2}\right)+\frac{x}{2}\mathbb{E}\left(\frac{1}{1-s(x)\beta_2}+\frac{1}{1+s(x)\beta_2}\right).
\end{equation}
On the other hand, we have for any $0\leq r\leq 1$
\[
\frac{1}{2}\mathbb{E}\left(\frac{1}{1-r\beta_2}+\frac{1}{1+r\beta_2}\right)=\mathbb{E}\frac{1}{1-\left(r\beta_2\right)^2}\leq \mathbb{E}\frac{1}{1-\beta_2^2}=\log 4,
\]
where the last equality follows from (\ref{uno}). This, together with (\ref{tres}), (\ref{cinco}), and (\ref{siete}), shows claim (\ref{cuatro}).

Finally, observe that
\[
\left\vert \Delta^2_{h\varphi(x)}f(x)\right\vert\leq 4\|B_nf-f\|+\left\vert \Delta^2_{h\varphi(x)}B_nf(x)\right\vert .
\]
Therefore, the result follows by applying claim (\ref{cuatro}) to $g=B_nf$.
\hfill \break \rightline{$\square$}

Let $(\widetilde{U}_k)_{k\geq 1}$ be a sequence of independent identically distributed random variables having the uniform distribution on $[0,1]$. Define
\begin{equation}\label{treceestrella}
S_n(x)=\sum_{k=1}^{n}1_{[0,x]}(\widetilde{U}_k),
\end{equation}
where $1_A$ stands for the indicator function of the set $A$. Clearly,
\begin{equation}\label{ocho}
P\left(S_n(x)=k\right)=\binom{n}{k}x^k(1-x)^{n-k},\quad k=0,1,\ldots n.
\end{equation}
We consider the sequence $L_n$ of positive linear operators  which associates to each function  $\phi:[0,\infty )\rightarrow \mathbb{R}$ the function
\begin{equation}\label{diecisiete}
L_n\phi(x)=\mathbb{E}\phi(S_n(x))=\sum_{k=0}^{n}\phi(k)P(S_n(x)=k).
\end{equation}
To obtain the derivatives of $L_n\phi$ in closed form, consider first the $m$th forward differences of $\phi$ at step $h\geq 0$, i.e.,
\begin{equation}\label{dieciocho}
\underline{\Delta}_h^m\phi(y)=\sum_{j=0}^m\binom{m}{j}(-1)^{m-j}\phi (y+hj),\quad y\geq 0,\quad m\in \mathbb{N}_0.
\end{equation}
If $\phi\in C^m[0,\infty )$, we have the representation (cf. \cite[Lemma 7.2]{adelllekuona})
\begin{equation}\label{diez}
\underline{\Delta}_h^m\phi(y)=h^m\mathbb{E}\phi^{(m)}\left(y+h\left(U_1+\cdots +U_m\right)\right),\quad y\geq 0,\quad n\in \mathbb{N},
\end{equation}
where $(U_k)_{k\geq 1}$ is a sequence of independent copies of a random variable $U$ having the uniform distribution on $[0,1]$.

In second place, denote by $K_m(x;y):=K_m(n,x;y)$, $m=0,1,\ldots,n$, the orthogonal polynomials with respect to the probability measure defined in (\ref{ocho}), namely, the  Krawtchouk polynomials. Such polynomials are explicitly defined (see, for instance, Chihara \cite[p. 161]{chihara})) by
\begin{equation}\label{diecinueve}
K_m(x;y)=\sum_{j=0}^m\binom{n-y}{m-j}\binom{y}{j}(-x)^{m-j}(1-x)^j,\quad m=0,1,\ldots,n,\quad y\in \mathbb{R},
\end{equation}
and satisfy the orthogonality property
\begin{equation}\label{once}
\mathbb{E}K_r(x;S_n(x))K_m(x;S_n(x))=\binom{n}{m}\varphi^{2m}(x)\delta_{r,m},\quad r,m\in \{0,1,\ldots n\},
\end{equation}
where $\delta_{r,m}$ is the Kronecker delta. In particular,
\begin{equation}\label{onceasterisco}
K_2(x;y)=\frac{1}{2}\left((y-nx)^2-(1-2x)(y-nx)-nx(1-x)\right).
\end{equation}

It was shown in \cite{adellanoz} (see also Roos \cite{roos} and L\'{o}pez-Bl\'{a}zquez and Salamanca \cite{lopezsalamanca}) that for $m=0,1,\ldots,n$
\begin{equation}\label{doce}
(L_n\phi)^{(m)}(x)=(n)_m\mathbb{E}\underline{\Delta}_1^m\phi(S_{n-m}(x))=\frac{m!}{\varphi^{2m}(x)}\mathbb{E}\phi(S_n(x))K_m(x;S_n(x)),
\end{equation}
where $(n)_m=n(n-1)\cdots (n-m+1)$.

As follows from (\ref{1}) and (\ref{ocho}), the $n$th Bernstein polynomial can be written in probabilistic terms as
\begin{equation}\label{trece}
B_nf(x)=\mathbb{E}f\left(\frac{S_n(x)}{n}\right).
\end{equation}
From now on, we assume that the random variables $S_n(x)$, $(U_k)_{k\geq 1}$, and $(\beta_m)_{m\geq 1}$, as defined in (\ref{ocho}), (\ref{diez}), and (\ref{uno}), respectively, are mutually independent.

\begin{proposition}\label{propo3}
Let $m\in \{1,2,\ldots,n\}$. Then,
\begin{equation}\label{catorce}
(B_nf)^{(m)}(x)=(n)_m\mathbb{E}\underline{\Delta}_{1/n}^mf\left(\frac{S_{n-m}(x)}{n}\right)=\frac{m!}{\varphi^{2m}(x)}\mathbb{E}f\left(\frac{S_n(x)}{n}\right)K_m(x;S_n(x)).
\end{equation}

In addition, if $f\in C^m[0,1]$, then
\begin{equation}\label{quince}
(B_nf)^{(m)}(x)=\frac{(n)_m}{n^m}\mathbb{E}f^{(m)}\left(\frac{S_{n-m}(x)+U_1+\cdots +U_m}{n}\right).
\end{equation}
\end{proposition}
\textit{Proof}.\
Set $\phi(y)=f(y/n)$, $0\leq y\leq n$. Observe that
\[
\underline{\Delta}_1^m\phi(y)=\underline{\Delta} _{1/n}^mf\left(\frac{y}{n}\right),\quad 0\leq y\leq n.
\]
Since $B_nf(x)=L_n\phi (x)$, this formula and (\ref{doce}) imply identity (\ref{catorce}). Formula (\ref{quince}) is an immediate consequence of (\ref{diez}) and the first equality in (\ref{catorce}). The proof is complete.
\hfill \break \rightline{$\square$}

We point out that the first equality in (\ref{catorce}) is well known (see, for instance, Abel et al. \cite{abelleviatanrasa}).  The right-hand side in (\ref{quince}) is the probabilistic representation of the Bernstein-Kantorovich operators acting on $f^{(m)}$ (see Acu et al. \cite{acurasasteopoaie} and the references therein).

\section{An auxiliary result for the central region}\label{sec3}

By $X \overset{\mathcal{L}}{=} Y$, we mean that the random variables $X$ and $Y$ have the same law. It follows from (\ref{ocho}) that
\begin{equation}\label{dieciseis}
S_n(1-x)\overset{\mathcal{L}}{=}n-S_n(x).
\end{equation}
On the other hand, to prove Theorem \ref{th1} it is enough to show that
$\left\|\varphi ^2(B_nf)^{(2)}\right\|\leq C_n\left\|B_nf-f\right\|$, $f\in \mathfrak{F}$, for some absolute constant $C_n$, as follows from Theorem \ref{th3}. We claim that it suffices to show that
\begin{equation}\label{diecisiete}
\sup_{0<x\leq 1/2} \left|\varphi ^2(x)(B_nf)^{(2)}(x)\right|\leq C_n \sup_{0<x\leq 1/2}\left|B_nf(x)-f(x)\right|,\quad f\in \mathfrak{F}.
\end{equation}
Actually, if $f\in \mathfrak{F}$, we define $\tilde{f}\in \mathfrak{F}$ as $\tilde{f}(y)=f(1-y)$, $y\in [0,1]$. By (\ref{onceasterisco}), the second equality in (\ref{catorce}), and (\ref{dieciseis}), it can be checked that
\[
B_nf(1-x)-f(1-x)=B_n\tilde{f}(x)-\tilde{f}(x),\qquad \varphi^2(1-x)(B_nf)^{(2)}(1-x)=\varphi^2(x)(B_n\tilde{f})^{(2)}(x),
\]
thus showing claim (\ref{diecisiete}). For this reason, we assume from now on that $0<x\leq 1/2$. For the sake of simplicity, we also denote by $\Vert \cdot \Vert$ the sup-norm on the interval $(0,1/2]$.

To show (\ref{diecisiete}), we distinguish if the supremum $\left\Vert \varphi^2(B_nf)^{(2)}\right\Vert$  is attained in the central region $[a_n,1/2]$ or in the non-central region $(0,a_n)$, for a suitable $0<a_n\leq 1/2$. The estimation procedures, inspired by the work of Sang\"{u}esa \cite{sanguesa}, are quite different in both regions.

From now on, we will use many times the inequality $1-y\leq e^{-y}$, $y\geq 0$, without explicitly mentioning it. The crucial quantity to estimate in the central region is the inverse moment
\begin{equation}\label{dieciocho}
H_n(x)=\frac{\varphi(x)\sqrt{n}}{n+2}\mathbb{E}\frac{|S_n(x)-nx|}{\varphi^2\left(\frac{S_n(x)+V}{n+2}\right)},
\end{equation}
where
\begin{equation}\label{diecinueve}
V=U_1+U_2,
\end{equation}
and $U_1$ and $U_2$ are the random variables defined in (\ref{diez}). Recall that the strong law of large numbers and the central limit theorem for $S_n(x)$ respectively state that
\[
\frac{S_n(x)}{n}\longrightarrow x,\ \text{a.s.},\qquad \text{and}\qquad \frac{S_n(x)-nx}{\varphi(x)\sqrt{n}}\overset{\mathcal{L}}{\longrightarrow} Z,\ \text{as}\ n\rightarrow \infty ,
\]
where $Z$ is a standard normal random variable. Thus, we have for a fixed $x\in (0,1/2]$
\[
\lim_{n\rightarrow \infty}H_n(x)=\lim_{n\rightarrow \infty}\mathbb{E}\frac{|S_n(x)-nx|}{\varphi(x)\sqrt{n}}=\mathbb{E}|Z|=\sqrt{\frac{2}{\pi}}= 0.7978\ldots
\]
This means that any upper bound for $H_n(x)$ cannot be better than $\sqrt{2/\pi}$. In fact, $H_n(x)$ takes bigger values when  $\lambda=nx$, for moderate values of $\lambda$, that is, when $S_n(x)$ approaches to a Poisson random variable. In this regard, let $N_{\lambda}$ be a random variable having the Poisson distribution with mean $\lambda$, i.e.,
\[
P\left(N_{\lambda}=k\right)=e^{-\lambda}\frac{\lambda^k}{k!},\quad k\in \mathbb{N}_0,\quad \lambda \geq 0.
\]
For any $\lambda \geq 0$, denote by
\begin{equation}\label{veinte}
C(\lambda)=2\log \frac{27}{16}\nu (\lambda)+r(\lambda),
\end{equation}
where
\begin{equation}\label{veintiuno}
\nu (\lambda)=\sqrt{\lambda}\left(2P\left(N_{\lambda}=
\lceil {\lambda}\rceil \right)-P\left(N_{\lambda}=
0\right)\right)+\frac{1}{\sqrt{\lambda}}\left(2P\left(N_{\lambda}\leq
\lceil {\lambda}\rceil \right)-1-P\left(N_{\lambda}=
0\right)\right),
\end{equation}
$\lceil {\lambda}\rceil$ being the ceiling of $\lambda$, and
\begin{equation}\label{veintidos}
r(\lambda)=\left(\log 4-2\log \frac{27}{16}\right)\lambda^{3/2}e^{-\lambda}.
\end{equation}
Numerical computations, carried out with the software Mathematica 10.2,  show that
\begin{equation}\label{veintitres}
\sup _{\lambda\geq 0}\ C(\lambda)=0.9827 \ldots < 0.99.
\end{equation}
Finally, let $\lambda_0>0$ be such that
\begin{equation}\label{veinticuatro}
\sqrt{\frac{2}{\pi}}+\frac{1}{\sqrt{\lambda_0}}\leq c:=0.8.
\end{equation}

We state the main result of this section.
\begin{theorem}\label{th4}
Suppose that $n\geq 2\lambda_0$. Then,
\[
H_n:=\sup_{0<x\leq 1/2}\ H_n(x)\leq 0.99+2\frac{D(\lambda_0)}{n}\log \frac{27}{16}+\frac{n^{3/2}}{2^{n+1/2}},
\]
where $D(\lambda_0)$ is the explicit constant defined in (\ref{treintaysiete}), only depending upon $\lambda_0$.
\end{theorem}

To prove this result, some auxiliary lemmas will be needed.

\begin{lemma}\label{lema5}
Let $V$ be as in (\ref{diecinueve}) and let $h(y)=y\log y$, $y\geq 0$. Then,
\begin{equation}\label{veinticinco}
\mathbb{E}\frac{1}{y+V}=\frac{1}{y+1}\sum_{j=0}^{\infty}\mathbb{E}\left(\frac{V-1}{y+1}\right)^{2j}=\underline{\Delta}_1^2h(y).
\end{equation}
In particular,
\[
\mathbb{E}\frac{1}{V}=\log 4,\quad \mathbb{E}\frac{1}{1+V}=\log \frac{27}{16}.
\]
\end{lemma}
\textit{Proof}.\
By (\ref{diecinueve}), the random variable $V$ is symmetric around $1$ and therefore $\mathbb{E}(V-1)^{2j+1}=0$, $j\in \mathbb{N}_0$. Hence,
\[
\mathbb{E}\frac{1}{y+V}=\frac{1}{y+1}\mathbb{E}\frac{1}{1+(V-1)/(y+1)}=\frac{1}{y+1}\sum_{j=0}^{\infty}\mathbb{E}\left(\frac{V-1}{y+1}\right)^{2j}.
\]
On the other hand, let $y>0$. Since $h^{(2)}(y)=1/y$, we have from (\ref{diez})
\[
\mathbb{E}\frac{1}{y+V}=\underline{\Delta}_1^2h(y).
\]
Letting $y\rightarrow 0$, this formula also holds for $y=0$. Finally, the second statement in Lemma \ref{lema5} follows by choosing $y=0$ and $y=1$ in (\ref{veinticinco}).
\hfill \break \rightline{$\square$}

Using (\ref{seis}), we rewrite the function $H_n(x)$ defined in (\ref{dieciocho}) as
\[
H_n(x)=\varphi(x)\sqrt{n}\left(\mathbb{E}\frac{|S_n(x)-nx|}{S_n(x)+V}+\mathbb{E}\frac{|S_n(x)-nx|}{n-S_n(x)+2-V}\right)
\]
\begin{equation}\label{veintiseis}
=\varphi(x)\sqrt{n}\left(\mathbb{E}\frac{|S_n(x)-nx|}{S_n(x)+V}+\mathbb{E}\frac{|S_n(1-x)-n(1-x)|}{S_n(1-x)+V}\right),
\end{equation}
as follows from (\ref{dieciseis}), the fact that $2-V\overset{\mathcal{L}}{=}V$, and the independence between $S_n(x)$ and $V$.

Observe that $H_n(x)$ depends on inverse moments involving the random variables $S_n(x)$ and $V$. In the following result, we estimate $H_n(x)$ in terms of inverse moments only involving $S_n(x)$.

\begin{lemma}\label{lema6}
Let $H_n(x)$ be as in (\ref{veintiseis}). Then,
\[
H_n(x)\leq 2\log \frac{27}{16}\left(I_n(x)+I_{n}(1-x)\right)+\left(\log 4-2\log \frac{27}{16}\right)(nx)^{3/2}(1-x)^{n+1/2}+\frac{n^{3/2}}{2^{n+1/2}},
\]
where
\begin{equation}\label{veintisiete}
I_n(x)=\varphi (x)\sqrt{n}\ \mathbb{E}\frac{|S_n(x)-nx|}{S_n(x)+1}.
\end{equation}
\end{lemma}
\textit{Proof}.\
Let us bound the first term on the right-hand side in (\ref{veintiseis}). Replace $y$ by $S_n(x)$ in (\ref{veinticinco}) and use the independence between $S_n(x)$ and $V$ to obtain
\begin{equation}\label{veintiocho}
\varphi (x)\sqrt{n}\ \mathbb{E}\frac{|S_n(x)-nx|}{S_n(x)+V}=\varphi (x)\sqrt{n}\left(\mathbb{E}\frac{|S_n(x)-nx|}{S_n(x)+1}+\sum_{j=1}^{\infty}\mathbb{E}(V-1)^{2j}\mathbb{E}\frac{|S_n(x)-nx|}{\left(S_n(x)+1\right)^{2j+1}}\right).
\end{equation}
Let $j\geq 1$. By considering the events $\{S_n(x)=0\}$ and $\{S_n(x)\geq 1\}$, we get
\[
\mathbb{E}\frac{|S_n(x)-nx|}{\left(S_n(x)+1\right)^{2j+1}}\leq nxP(S_n(x)=0)+\frac{1}{2^{2j}}\mathbb{E}\frac{|S_n(x)-nx|}{S_n(x)+1}1_{\{S_n(x)\geq 1\}}
\]
\[
=nxP(S_n(x)=0)+\frac{1}{2^{2j}}\left(\mathbb{E}\frac{|S_n(x)-nx|}{S_n(x)+1}-nxP(S_n(x)=0)\right).
\]
This and Lemma \ref{lema5} imply that the sum in (\ref{veintiocho}) is bounded above by
\[
(\log 4-1)nxP(S_n(x)=0)+\left(2\log \frac{27}{16}-1\right)\left(\mathbb{E}\frac{|S_n(x)-nx|}{S_n(x)+1}-nxP(S_n(x)=0)\right).
\]
We therefore have from (\ref{veintisiete}) and (\ref{veintiocho})
\begin{equation}\label{veintinueve}
\varphi (x)\sqrt{n}\ \mathbb{E}\frac{|S_n(x)-nx|}{S_n(x)+V}\leq 2\log \frac{27}{16}I_n(x)+\left(\log 4-2\log \frac{27}{16}\right)(nx)^3(1-x)^{n+1/2}.
\end{equation}
To bound the second term on the right-hand side in (\ref{veintiseis}), we follow the same procedure, replacing $x$ by $1-x$ and noting that
\[
(n(1-x))^{3/2}x^{n+1/2}\leq \frac{n^{3/2}}{2^{n+1/2}},
\]
since, by assumption, $0<x\leq 1/2$. This and (\ref{veintinueve}) show the result.
\hfill \break \rightline{$\square$}

An exact expression for $I_n(x)$, useful to perform computations, is given in the following result.

\begin{lemma}\label{lema7}
Let $I_n(x)$ be as in (\ref{veintisiete}). Set $\lambda=nx$. Then,
\[
I_n(x)=\frac{n(1-x)^{3/2}}{n+1}\left(\sqrt{\lambda}\left(2P(S_n(x)=\lceil \lambda\rceil)-P(S_n(x)=0)\right)+\frac{1}{\sqrt{\lambda}}\left(2P(S_n(x)\leq \lceil \lambda \rceil )-1-P(S_n(x)=0)\right)\right).
\]
\end{lemma}
\textit{Proof}.\
Since $|y|=2y_+-y$, where $y_+=\max (0,y)$, $y\in \mathbb{R}$, we can write
\begin{equation}\label{treinta}
\mathbb{E}\frac{|S_n(x)-nx|}{S_n(x)+1}=2\mathbb{E}\frac{(S_n(x)-nx)_+}{S_n(x)+1}-\mathbb{E}\frac{S_n(x)-nx}{S_n(x)+1}.
\end{equation}
Using the identity
\begin{equation}\label{treintayuno}
\frac{1}{k+1}P(S_n(x)=k)=\frac{1}{(n+1)x}P(S_{n+1}(x)=k+1),\quad k=0,1,\ldots ,n,
\end{equation}
we have
\[
\mathbb{E}\frac{S_n(x)-nx}{S_n(x)+1}=\sum _{k=0}^n\frac{k-\lambda}{k+1}P(S_n(x)=k)=1-\frac{\lambda+1}{(n+1)x}P(S_{n+1}(x)\geq 1)
\]
\begin{equation}\label{treintaydos}
=-\frac{1-x}{(n+1)x}+\frac{\lambda+1}{(n+1)x}P(S_{n+1}(x)=0)=-\frac{1-x}{(n+1)x}\left(1-(\lambda+1)P(S_n(x)=0)\right).
\end{equation}
Also, by (\ref{treceestrella}) and (\ref{treintayuno}), we get
\[
\mathbb{E}\frac{(S_n(x)-nx)_+}{S_n(x)+1}=\sum _{k=\lceil\lambda \rceil}^n\frac{k-\lambda}{k+1}P(S_n(x)=k)=P(S_n(x)\geq \lceil \lambda \rceil)-\frac{\lambda+1}{(n+1)x}P(S_{n+1}(x)\geq \lceil \lambda \rceil +1)
\]
\[
=P(S_n(x)\geq \lceil \lambda \rceil )-\frac{\lambda+1}{(n+1)x}\left(xP(S_n(x)\geq \lceil \lambda \rceil)+(1-x)P(S_n(x)\geq \lceil \lambda \rceil +1)\right)
\]
\[
=\frac{1-x}{(n+1)x}\left(\lambda P(S_n(x)\geq \lceil \lambda \rceil)-(\lambda +1)P(S_n(x)\geq \lceil \lambda \rceil +1)\right)
\]
\[
=\frac{1-x}{(n+1)x}\left(\lambda P(S_n(x)= \lceil \lambda \rceil)-P(S_n(x)\geq \lceil \lambda \rceil +1)\right).
\]
This, together with (\ref{treinta}), (\ref{treintaydos}), and some simple computations, shows the result.
\hfill \break \rightline{$\square$}

We give two upper  bounds for $I_n(x)$, according to whether or not $\lambda=nx$ is large. In the first case, we use the well known Stirling's approximation stating that
\[
e^{1/(2m+1)}\sqrt{2\pi m}\leq \frac{m!}{m^m}e^m\leq e^{1/(2m)}\sqrt{2\pi m},\quad m\in \mathbb{N}.
\]
Observe that this implies that
\begin{equation}\label{treintaytres}
P(S_n(m/n)=m)\leq \frac{1}{\sqrt{2\pi}}\sqrt{\frac{n}{m(n-m)}},\quad m=1,\ldots, n-1.
\end{equation}
For the second case, recall that the total variation distance between two $\mathbb{N}_0$-valued random variables $X$ and $Y$ is defined as
\begin{equation}\label{cuarentaydosestrella}
d_{TV}(X,Y)=\sup _{A\subseteq \mathbb{N}_0}\ |P(X\in A)-P(Y\in A)|=\frac{1}{2}\sum_{k=0}^{\infty }|P(X=k)-P(Y=k)|.
\end{equation}
It was shown in \cite[formulae (1.7) and (2.20)]{adellanozlekuona} that
\begin{equation}\label{treintaycuatro}
d_{TV}(S_n(x),N_{\lambda})\leq \frac{\lambda}{n}\left(\frac{\sqrt{2}}{  4}+\frac{4}{11}(3\lambda+4)\frac{\lambda ^2}{n}\right),\quad \lambda=nx,\quad n\geq 10.
\end{equation}
Other estimates can be found in Barbour and Hall \cite{barbourhall} and Deheuvels et al. \cite{deheuvelspfeifer}.

\begin{lemma}\label{lema8}
Let $\lambda_0$ be as in (\ref{veinticuatro}) and assume that $n\geq 2\lambda_0$. If $\lambda:=nx\geq \lambda_0$, then
\begin{equation}\label{treintaycinco}
I_n(x)\leq c(1-x),
\end{equation}
where $c$ is defined in (\ref{veinticuatro}). If $\lambda =nx < \lambda_0$, then
\begin{equation}\label{treintayseis}
I_n(x)\leq (1-x)\left(\nu (\lambda)+\frac{D(\lambda_0)}{n}\right),
\end{equation}
where $\nu (\lambda)$ is defined in (\ref{veintiuno}) and
\begin{equation}\label{treintaysiete}
D(\lambda_0)=3\sqrt{\lambda_0}(\lambda_0+1)\left(\frac{\sqrt{2}}{4}+\frac{2}{11}(3\lambda_0+4)\lambda_0\right).
\end{equation}
\end{lemma}
\textit{Proof}.\  Assume that $\lambda \geq \lambda_0$. From Lemma \ref{lema7} and (\ref{veinticuatro}), we see that
\begin{equation}\label{treintayocho}
I_n(x)\leq (1-x)\left(2\sqrt{\frac{\lambda(n-\lambda)}{n}}P((S_n(x)=\lceil \lambda \rceil )+\frac{1}{\sqrt{\lambda_0}}\right).
\end{equation}
Differentiating with respect to $x$, it can be checked that
\[
P((S_n(x)=\lceil \lambda \rceil )\leq P((S_n(\lceil \lambda \rceil /n)=\lceil \lambda \rceil ).
\]
We thus have from (\ref{treintaytres})
\begin{equation}\label{treintaynueve}
2\sqrt{\frac{\lambda(n-\lambda)}{n}}P(S_n(x)=\lceil \lambda \rceil )\leq \sqrt{\frac{2}{\pi}}\sqrt{\frac{\lambda(n-\lambda)}{\lceil \lambda \rceil (n-\lceil \lambda \rceil )}}\leq \sqrt{\frac{2}{\pi}},
\end{equation}
where the last inequality follows from two facts: first $n/2\geq \lambda \geq \lambda_0$, since $0<x\leq 1/2$, and second, the function $h(\lambda)=\lambda (n-\lambda)$ increases for $n\geq 2\lambda$. Therefore, (\ref{treintaycinco}) follows from (\ref{veinticuatro}), (\ref{treintayocho}), and (\ref{treintaynueve}).

Since $\lambda_0>5$, inequality (\ref{treintayseis}) is an immediate consequence of Lemma \ref{lema7} and (\ref{treintaycuatro}).
\hfill \break \rightline{$\square$}

\textit{Proof of Theorem \ref{th4}}
Let $C(\lambda)$ be as in (\ref{veinte}) and denote
\[
\widetilde{C}(\lambda)=2c\log \frac{27}{16}+r(\lambda),\quad \lambda \geq 0,
\]
where $c=0.8$ and $r(\lambda)$ is defined in (\ref{veintidos}). Numerical computations show that
\begin{equation}\label{cuarenta}
0.9792\ldots =\sup _{\lambda\geq 0}\widetilde{C}(\lambda)<\sup _{\lambda\geq 0}C(\lambda)=0.9827\ldots <0.99.
\end{equation}
Since $n(1-x)\geq n/2\geq \lambda_0$, we have from (\ref{treintaycinco})
\begin{equation}\label{cuarentayuno}
I_n(1-x)\leq cx.
\end{equation}
Set $\lambda =nx$. We distinguish the following two cases.

\underline{Case $\lambda\geq \lambda_0$}. We have from (\ref{treintaycinco}), (\ref{cuarentayuno}), and Lemma \ref{lema6}
\begin{equation}\label{cuarentaydos}
H_n(x)\leq 2c\log \frac{27}{16}+r(\lambda)+\frac{n^{3/2}}{2^{n+1/2}}=\widetilde{C}(\lambda )+\frac{n^{3/2}}{2^{n+1/2}}.
\end{equation}
\underline{Case $\lambda < \lambda_0$}. We have from (\ref{treintayseis}), (\ref{cuarentayuno}), and Lemma \ref{lema6}
\[
H_n(x)\leq 2\log \frac{27}{16}\left(\nu (\lambda)(1-x)+cx\right)+r(\lambda)+2\log \frac{27}{16}\frac{D(\lambda_0)}{n}+\frac{n^{3/2}}{2^{n+1/2}}
\]
\begin{equation}\label{cuarentaytres}
\leq \max (C(\lambda),\widetilde{C}(\lambda))+2\log \frac{27}{16}\frac{D(\lambda_0)}{n}+\frac{n^{3/2}}{2^{n+1/2}}.
\end{equation}
In view of (\ref{cuarenta}), the conclusion follows from (\ref{cuarentaydos}) and (\ref{cuarentaytres}).
\hfill \break \rightline{$\square$}

\section{Estimates in the central region}\label{sec4}

We define the central region as
\[
R_n(a):=\{x\in (0,1/2]:\ n\varphi ^2(x)>a \},
\]
where $a>0$ will be chosen later on. We make the following assumption

($H_1$)\qquad \qquad
$\displaystyle{
\left\Vert \varphi ^2(B_nf)^{(2)}\right\Vert =\sup \left\{ \left\vert \varphi ^2(x)(B_nf)^{(2)}(x)\right\vert :\ x\in R_n(a)\right\}.}
$

We denote by $B_n^k$, $k\in \mathbb{N}$, the $k$th iterate of the operator $B_n$, as well as
\begin{equation}\label{cuarentaycuatro}
    g=B_nf,\qquad \tau_m(x)=x+\left(\frac{S_n(x)}{n}-x\right)\beta_m,\quad m\in \mathbb{N}.
\end{equation}
Since
\begin{equation}\label{cuarentaycinco}
   \mathbb{E}\left(\frac{S_n(x)}{n}-x\right)^2=\frac{\varphi^2(x)}{n},
\end{equation}
we use Taylor's formula (\ref{dos}) replacing $y$ by $S_n(x)/n$ to obtain the following starting identity for the central region
\[
B_n^3f(x)-B_n^2f(x)=B_n^2g(x)-B_ng(x)=\mathbb{E}B_ng\left(\frac{S_n(x)}{n}\right)-B_ng(x)
\]
\begin{equation}\label{cuarentayseis}
   =\frac{\varphi^2(x)}{2n}g^{(2)}(x)+\frac{\varphi^2(x)}{2n}\left((B_ng)^{(2)}(x)-g^{(2)}(x)\right)+\frac{1}{6}\mathbb{E}(B_ng)^{(3)}(\tau_3(x))\left(\frac{S_n(x)}{n}-x\right)^{3}.
\end{equation}
The second term on the right-hand side in (\ref{cuarentayseis}) is bounded as follows.

\begin{lemma}\label{lema9}
We have
\[
\frac{1}{2n}\left\Vert \varphi^2\left((B_ng)^{(2)}-g^{(2)}\right)\right\Vert \leq \frac{1}{\sqrt{2}}\left\Vert B_nf-f\right \Vert .
\]
\end{lemma}
\textit{Proof}.\  By (\ref{cuarentaycuatro}) and the second equality in (\ref{catorce}), we have
\[
\frac{\varphi ^2(x)}{2n}\left\vert (B_ng)^{(2)}(x)-g^{(2)}(x) \right\vert = \frac{\varphi ^2(x)}{2n}\left\vert \left( B_n^2f-B_nf
 \right)^{(2)}(x) \right\vert
\]
\[
=\frac{1}{n\varphi^2(x)}\left\vert \mathbb{E}(B_nf-f)\left(\frac{S_n(x)}{n}\right)K_2(x;S_n(x))\right\vert \leq  \frac{1}{n\varphi ^2(x)} \Vert B_nf-f\Vert \mathbb{E}\left\vert K_2(x;S_n(x)) \right \vert .
\]
Therefore, the result follows from (\ref{once}) and Schwarz's inequality.
\hfill \break \rightline{$\square$}

To estimate the last term in (\ref{cuarentayseis}), some auxiliary results will be needed.

\begin{lemma}\label{lema10}
Let $H_n$ be as in Theorem \ref{th4}. Then,
\[
\left\Vert \varphi^3 (B_ng)^{(3)} \right\Vert \leq \sqrt{n+1}H_{n-2}\left\Vert \varphi ^2g^{(2)}\right\Vert .
\]
\end{lemma}
\textit{Proof}.\  Applying (\ref{quince}) with $m=2$, we have
\begin{equation}\label{cuarentaysiete}
  (B_ng)^{(2)}(x)=\frac{n-1}{n}\mathbb{E}g^{(2)}\left(\frac{S_{n-2}(x)+V}{n}\right),
\end{equation}
where $V$ is defined in (\ref{diecinueve}). We differentiate this formula using the second equality in (\ref{catorce}) and take into account that $K_1(x;y)=y-nx$ to obtain
\[
(B_ng)^{(3)}(x)=\frac{n-1}{n} \frac{1}{\varphi ^2(x)}\mathbb{E}g^{(2)}\left(\frac{S_{n-2}(x)+V}{n}\right)(S_{n-2}(x)-(n-2)x).
\]
By (\ref{dieciocho}), this implies that
\[
\left\vert \varphi ^3(x)(B_ng)^{(3)}(x)\right\vert \leq \frac{n-1}{n}\left\Vert \varphi ^2g^{(2)}\right\Vert \varphi (x)\mathbb{E}\left(\frac{\left\vert S_{n-2}(x)-(n-2)x\right\vert}{\varphi ^2\left(\frac{S_{n-2}(x)+V}{n}\right)}\right)=\frac{n-1}{\sqrt{n-2}}\left\Vert \varphi ^2g^{(2)}\right\Vert H_{n-2}(x),
\]
which shows the result, because $n\geq 3$.
\hfill \break \rightline{$\square$}

Denote by
\[
\mu_n(x)=\mathbb{E}\left\vert \frac{S_n(x)}{n}-x\right\vert ^k,\quad k\in \mathbb{N},
\]
the $k$th absolute central moment of the Bernstein polynomials. Explicit expressions for the first even central moments can be found, for instance, in \cite{rm2019}. Particularly,
\begin{equation}\label{cuarentayocho}
 \mu _4(x)=\frac{3\varphi^4(x)}{n^2}+\frac{\varphi^2(x)}{n^3}(1-6\varphi^2(x))\leq \frac{\varphi^4(x)}{n^2}\left(3+\frac{1}{n\varphi^2(x)}\right),
\end{equation}
and
\[
 \mu _6(x)=\frac{\varphi^2(x)}{n^5}\left(15\varphi^4(x)n^2+5\varphi^2(x)\left(5-26\varphi^2(x)\right)n+1-30\varphi^2(x)(1-2x)^2\right)
\]
\begin{equation}\label{cuarentaynueve}
\leq \frac{\varphi^6(x)}{n^3}\left( 15+\frac{25}{n\varphi^2(x)}+\frac{1}{(n\varphi^2(x))^2}\right) .
\end{equation}

\begin{lemma}\label{lema11}
We have
\[
\frac{n\sqrt{n}}{\varphi^3(x)}\left(\mu _3(x)+\frac{3}{8}\frac{\mu_4(x)}{\varphi^2(x)}+\frac{3}{16}\frac{\mu_5(x)}{\varphi^4(x)} +\frac{7}{16}\frac{\mu_6(x)}{\varphi^6(x)}\right) \leq K(n\varphi^2(x)),
\]
where
\[
K(s)=\sqrt{3+\frac{1}{s}}+\frac{3}{8\sqrt{s}}\left(3+\frac{1}{s}\right)+\frac{3}{16s}\sqrt{\left(3+\frac{1}{s}\right)\left(15+\frac{25}{s}+\frac{1}{s^2}\right)}
\]
\begin{equation}\label{cincuenta}
+\frac{7}{16s\sqrt{s}}\left(15+\frac{25}{s}+\frac{1}{s^2}\right),\quad s>0.
\end{equation}
\end{lemma}
\textit{Proof}.\  By Schwarz's inequality, we see that
\[
\mu _3(x)=\mathbb{E}\left\vert \frac{S_n(x)}{n}-x\right\vert \left\vert \frac{S_n(x)}{n}-x\right\vert ^2 \leq \sqrt{\mu_2(x)\mu _4(x)}.
\]
Similarly,
\[
\mu_5(x)\leq  \sqrt{\mu_4(x)\mu _6(x)}.
\]
Thus, the result follows from (\ref{cuarentaycinco}), (\ref{cuarentayocho}), (\ref{cuarentaynueve}), and some simple computations.
\hfill \break \rightline{$\square$}

\begin{lemma}\label{lema12}
Let $\beta _m$ be as in (\ref{uno}), $m\in \mathbb{N}$. For any $0\leq z\leq 1$, we have
\[
\mathbb{E}\frac{\varphi^m(x)}{\varphi^m(x+(z-x)\beta_m)}\leq 1+\frac{m}{2(m+1)}\frac{|z-x|}{\varphi^2(x)}+\frac{m}{4(m+1)}\frac{|z-x|^2}{\varphi^4(x)}+\frac{m+4}{4(m+1)}\frac{|z-x|^3}{\varphi^6(x)}.
\]
\end{lemma}
\textit{Proof}.\  We claim that
\begin{equation}\label{cincuentayuno}
\mathbb{E}\frac{1}{(1-r\beta_m)^{m/2}}
\leq 1+\frac{m}{2(m+1)}r+\frac{m}{4(m+1)}r^2+\frac{m+4}{4(m+1)}r^3,\quad 0\leq r\leq 1.
\end{equation}
Indeed, we see from (\ref{uno}) that
\begin{equation}\label{cincuentaydos}
\mathbb{E}\frac{1}{(1-\beta_m)^{m/2}}=m\int _0^1(1-\theta)^{m/2-1}d\theta =2.
\end{equation}
Using the binomial expansion and the fact that $\mathbb{E}\beta_m=1/(m+1)$ and $\mathbb{E}\beta_m^2=2/((m+1)(m+2))$, we have
\[
\mathbb{E}\frac{1}{(1-r\beta_m)^{m/2}}
\leq 1+\frac{m}{2(m+1)}r+\frac{m}{4(m+1)}r^2+r^3\sum_{k=3}^{\infty }(-1)^k{-m/2 \choose k}\mathbb{E}\beta_m^k
\]
\begin{equation}\label{cincuentaytres}
=1+\frac{m}{2(m+1)}r+\frac{m}{4(m+1)}r^2+r^3\left(-1-\frac{m}{2(m+1)}-\frac{m}{4(m+1)}+\sum_{k=0}^{\infty }(-1)^k{-m/2 \choose k}\mathbb{E}\beta_m^k\right).
\end{equation}
This shows claim (\ref{cincuentayuno}), since the series in (\ref{cincuentaytres}) equals $2$, as follows from (\ref{cincuentaydos}).

Assume first that $0\leq z\leq x$. As noted in Sang\"{u}esa \cite{sanguesa}, a useful property of the weight function $\varphi$ is
\[
\varphi^2(rx)\geq r\varphi^2(x),\quad 0\leq r\leq 1.
\]
As a consequence,
\[
\varphi^m(x+(z-x)\beta_m)=\varphi^m\left(x\left(1-\frac{x-z}{x}\beta_m\right)\right)\geq \left(1-\frac{x-z}{x}\beta_m\right)^{m/2}\varphi^m(x).
\]
We therefore have from (\ref{cincuentayuno})
\[
\mathbb{E}\frac{\varphi^m(x)}{\varphi^m(x+(z-x)\beta_m)}\leq \mathbb{E}\frac{1}{\left(1-\frac{x-z}{x}\beta_m\right)^{m/2}}
\]
\[
\leq 1+\frac{m}{2(m+1)}\frac{|z-x|}{\varphi^2(x)}(1-x)+\frac{m}{4(m+1)}\frac{|z-x|^2}{\varphi^4(x)}(1-x)^2+\frac{m+4}{4(m+1)}\frac{|z-x|^3}{\varphi^6(x)}(1-x)^3,
\]
thus showing the result in this case. If $x\leq z\leq 1$, the proof is similar by observing that $\varphi^2(x)=\varphi^2(1-x)$ and
\[
\varphi^2(x+(z-x)\beta_m)=\varphi^2(1-x-(z-x)\beta_m)=\varphi^2\left((1-x)\left(1-\frac{z-x}{1-x}\beta_m\right)\right).
\]
This completes the proof.
\hfill \break \rightline{$\square$}

We are in a position to estimate the last term in (\ref{cuarentayseis}).

\begin{lemma}\label{lema13}
Let $H_n$ and $K(s)$ be as in Theorem \ref{th4} and (\ref{cincuenta}), respectively. Then,
\[
\frac{1}{6}\left\vert \mathbb{E}(B_ng)^{(3)}(\tau_3(x))\left(\frac{S_n(x)}{n}-x\right)^3\right\vert \leq \sqrt{\frac{n+1}{n}}\frac{H_{n-2}K(n\varphi^2(x))}{3}\frac{\left\Vert \varphi^2g^{(2)}\right\Vert}{2n}.
\]
\end{lemma}
\textit{Proof}.\  Fix $z=S_n(x)/n$. Applying Lemma \ref{lema10} and Lemma \ref{lema12} with $m=3$, we get
\[
\left\vert \mathbb{E}(B_ng)^{(3)}(x+(z-x)\beta_3)\right\vert \leq \frac{\left\Vert \varphi^3(B_ng)^{(3)}\right\Vert}{\varphi^3(x)}\mathbb{E}\frac{\varphi^3(x)}{\varphi^3(x+(z-x)\beta _3)}
\]
\[
\leq  \frac{\sqrt{n+1}H_{n-2}\left\Vert \varphi^2g^{(2)}\right\Vert }{\varphi^3(x)}\left(1+\frac{3}{8}\frac{|z-x|}{\varphi^2(x)}+\frac{3}{16}\frac{|z-x|^2}{\varphi^4(x)}+\frac{7}{16}\frac{|z-x|^3}{\varphi^6(x)}\right).
\]
We multiply both sides of this inequality by $|z-x|^3/6$ and replace $z$ by $S_n(x)/n$. Recalling (\ref{cuarentaycuatro}), we obtain
\[
\frac{1}{6}\left\vert \mathbb{E}(B_ng)^{(3)}(\tau_3(x))\left(\frac{S_n(x)}{n}-x\right)^3\right\vert
\]
\[
 \leq \sqrt{\frac{n+1}{n}}\frac{H_{n-2}}{3}\frac{\left\Vert \varphi^2g^{(2)}\right\Vert}{2n}\frac{n\sqrt{n}}{\varphi^3(x)}\left(\mu_3(x)+\frac{3}{8}\frac{\mu_4(x)}{\varphi^2(x)}+\frac{3}{16}\frac{\mu_5(x)}{\varphi^4(x)}+\frac{7}{16}\frac{\mu_6(x)}{\varphi^6(x)}\right).
\]
Thus, the conclusion follows from Lemma \ref{lema11}.
\hfill \break \rightline{$\square$}

We state the main result of this section.

\begin{theorem}\label{th14}
Assume ($H_1$). Then,
\[
\left(1-\sqrt{\frac{n+1}{n}}\frac{H_{n-2}K(a)}{3}\right)\frac{\left\Vert \varphi^2g^{(2)}\right\Vert}{2n}\leq \frac{\sqrt{2}+1 }{\sqrt{2}}\Vert B_nf-f \Vert .
\]
\end{theorem}
\textit{Proof}.\  Let $x\in R_n(a)$. Starting from (\ref{cuarentayseis}) and applying Lemma \ref{lema9}, we obtain
\[
\frac{\left\vert \varphi^2(x)g^{(2)}(x)\right\vert}{2n}\leq \Vert B_nf-f\Vert+\frac{1}{\sqrt{2}}\Vert B_nf-f\Vert +\frac{1}{6}\left\vert \mathbb{E}(B_ng)^{(3)}(\tau_3(x))\left(\frac{S_n(x)}{n}-x\right)^3\right\vert
\]
\[
\leq \frac{\sqrt{2}+1}{\sqrt{2}}\Vert B_nf-f\Vert +\sqrt{\frac{n+1}{n}}\frac{H_{n-2}K(a)}{3}\frac{\left\Vert \varphi^2g^{(2)}\right\Vert}{2n},
\]
where the last inequality follows from Lemma \ref{lema13}. Thus, the result follows from assumption ($H_1$).
\hfill \break \rightline{$\square$}

Theorem \ref{th14} makes sense only when the term in front of $\left\Vert \varphi^2g^{(2)}\right\Vert/(2n)$  is positive. This, in turn, implies that the parameter $a$ cannot be close to $0$, because of the definition of $K(a)$ given in (\ref{cincuenta}).

\section{Estimates in the non-central region}\label{sec5}

We define the non-central region as
\[
\overline{R}_n(a):=\{x\in (0,1/2]:\ n\varphi ^2(x)<a \},\quad a>0,
\]
and make the alternative assumption

($H_2$)\qquad \qquad
$\displaystyle{
\left\Vert \varphi ^2g^{(2)}\right\Vert=\left\Vert \varphi ^2(B_nf)^{(2)}\right\Vert =\sup \left\{ \left\vert \varphi ^2(x)(B_nf)^{(2)}(x)\right\vert :\ x\in \overline{R}_n(a)\right\}.}
$

In this section, we assume that $x\in \overline{R}_n(a)$. The starting point is the following identity involving the second derivatives of the iterates of the Bernstein polynomials
\begin{equation}\label{cincuentaycuatro}
\varphi^2(x)g^{(2)}(x)=\sum_{k=0}^m\varphi^2(x)\left( (B_n^{k+1}g)^{(2)}(x)-(B_n^kg)^{(2)} (x)\right)+\varphi^2(x)(B_n^{m+1}g)^{(2)}(x),\quad m\in \mathbb{N}.
\end{equation}

We recall here the following probabilistic representation of the terms in (\ref{cincuentaycuatro}) given in Sang\"{u}esa \cite{sanguesa}. Denote
\begin{equation}\label{cincuentaycinco}
W_n(\theta)=\frac{S_{n-2}(\theta)+V}{n},\quad 0\leq \theta \leq 1,
\end{equation}
where $V$ is defined in (\ref{diecinueve}). Let $(S_{n-2,i}(\theta),\ 0\leq \theta \leq 1)_{i\geq 1}$ and $(V_i)_{i\geq 1}$ be sequences of independent copies of $(S_{n-2}(\theta),\ 0\leq \theta \leq 1)$ and $V$, respectively. Define
\[
W_{n,i}(\theta)=\frac{S_{n-2,i}(\theta)+V_i}{n},\quad 0\leq \theta \leq 1,\quad i\in \mathbb{N}.
\]
Clearly, $(W_{n,i}(\theta),\ 0\leq \theta \leq 1)_{i\geq 1}$ is a sequence of independent copies of $(W_n(\theta),\ 0\leq \theta \leq 1)$, as defined in (\ref{cincuentaycinco}). We consider the subordinated stochastic processes $(W_n^{(m)}(\theta),\ 0\leq \theta \leq 1)_{m\geq 1}$ inductively defined as follow
\begin{equation}\label{cincuentayseis}
W_n^{(1)}(\theta)=W_{n,1}(\theta),\qquad W_n^{(m+1)}(\theta)=W_{n,m+1}\left(W_n^{(m)}(\theta)\right),\quad 0\leq \theta \leq 1,\quad m\in \mathbb{N}.
\end{equation}

The following result was shown  by Sang\"{u}esa \cite{sanguesa}. We give here a short proof of it for the sake of completeness.

\begin{lemma}\label{lema15}
Let $m\in \mathbb{N}_0$. Then,
\[
(B_n^{m+1}g)^{(2)}(x)=\left(\frac{n-1}{n}\right)^{m+1} \mathbb{E}g^{(2)}\left(W_n^{(m+1)}(x)\right).
\]
\end{lemma}
\textit{Proof}.\  For $m=0$, the result follows from (\ref{cuarentaysiete}) and (\ref{cincuentayseis}). Let $0\leq \theta \leq 1$. Suppose that
\[
(B_n^mg)^{(2)}(\theta)=\left(\frac{n-1}{n}\right)^m\mathbb{E}g^{(2)}\left(W_n^{(m)}(\theta)\right),
\]
for some $m\in \mathbb{N}$. Again by (\ref{cuarentaysiete}) and (\ref{cincuentayseis}), we have
\[
(B_n^{m+1}g)^{(2)}(\theta)=\left(\frac{n-1}{n}\right)^m\mathbb{E}(B_ng)^{(2)}\left(W_n^{(m)}(\theta)\right)=\left(\frac{n-1}{n}\right)^{m+1}\mathbb{E}g^{(2)}\left(  W_{n,m+1}\left(W_n^{(m)}(\theta)\right)   \right),
\]
\[
=\left(\frac{n-1}{n}\right)^{m+1}\mathbb{E}g^{(2)}\left(  W_n^{(m+1)}(\theta)  \right),
\]
thus completing the proof.
\hfill \break \rightline{$\square$}

Observe that if $x\in \overline{R}_n(a)$, then
\begin{equation}\label{cincuentaysiete}
nx<b_n:=\frac{2a}{1+\sqrt{1-4a/n}},\quad n>4a.
\end{equation}
For this reason, we assume in this section that $n>4a$. The main quantity to be estimated is the inverse moment.
\begin{equation}\label{cincuentayocho}
\frac{1}{n}\mathbb{E}\frac{1}{\varphi^2\left(W_n^{(m)}(x)\right)},\quad m\in \mathbb{N}.
\end{equation}
To bound (\ref{cincuentayocho}), we will proceed by induction on $m$. In this respect, we define the iterates
\begin{equation}\label{sesentaytres}
\alpha_0(\theta)=\theta,\quad \alpha_1(\theta)=1-e^{-\theta},\quad \alpha_{m+1}(\theta)=\alpha_1(\alpha_m(\theta)),\quad 0\leq \theta\leq 1,\quad m\in \mathbb{N}.
\end{equation}
Note that $\alpha_m(\theta)\rightarrow 0$, as $m\rightarrow \infty $, for any $0<\theta \leq 1$. However, the speed of convergence is very slow, since $\alpha_1^{(1)}(0)=1$. The reason for these iterates comes from the following inequality, which follows from (\ref{cincuentaysiete}):
\begin{equation}\label{sesentaycuatro}
\mathbb{E}e^{-\theta S_{n-2}(x)}=\left(1+x(1-e^{-\theta})\right)^{n-2}\leq e^{-(n-2)x\alpha_1(\theta)}\leq e^{2b_n/n}e^{-nx\alpha_1(\theta)},\quad 0\leq \theta\leq 1.
\end{equation}
Also, we have from (\ref{diecinueve})
\begin{equation}\label{sesentaycinco}
\mathbb{E}e^{-\theta V}=\left(\mathbb{E}e^{-\theta U_1}\right)^2=\left(\frac{1-e^{-\theta}}{\theta}\right)^2=\left(\frac{\alpha_1(\theta)}{\theta}\right)^2,\quad 0\leq \theta\leq 1.
\end{equation}

\begin{lemma}\label{lema17}
For any $m\in \mathbb{N}$, we have
\begin{align*}
& \frac{1}{n}\mathbb{E}\frac{1}{\varphi^2\left(W_n^{(m)}(x)\right)}
\cr
& \leq e^{2b_n/n}\left(2\log \frac{27}{16}\int_0^1\left(\frac{\alpha_{m-1}(\theta)}{\theta}\right)^2e^{-nx\alpha_{m-1}(\theta)}d\theta
+\left(\log 4-2\log \frac{27}{16}\right)\alpha_{m-1}^2(1)e^{-nx\alpha_{m-1}(1)}+\epsilon _n\right),
\end{align*}
where $b_n$ is defined in (\ref{cincuentaysiete}) and
\begin{equation}\label{cincuentaynueve}
\epsilon_n=\frac{4}{n}\log \frac{27}{16}+e^{-n/2}.
\end{equation}
\end{lemma}

\textit{Proof}.\  We start with the case $m=1$. Using (\ref{seis}), we have as in (\ref{veintiseis})
\[
\frac{1}{n}\mathbb{E}\frac{1}{\varphi^2\left(W_n^{(1)}(x)\right)}=\frac{1}{n}\mathbb{E}\frac{1}{\varphi^2\left(\frac{S_{n-2}(x)+V}{n}\right)}=\mathbb{E}\frac{1}{S_{n-2}(x)+V}+\mathbb{E}\frac{1}{n-2-S_{n-2}(x)+2-V}
\]
\begin{equation}\label{sesenta}
=\mathbb{E}\frac{1}{S_{n-2}(x)+V}+\mathbb{E}\frac{1}{S_{n-2}(1-x)+V}.
\end{equation}
We claim that
\begin{equation}\label{sesentayuno}
\mathbb{E}\frac{1}{S_{n-2}(x)+V}\leq e^{2b_n/n}\left(2\log \frac{27}{16}\int_0^1e^{-nx\theta}d\theta +\left(\log 4-2\log \frac{27}{16}\right)e^{-nx}\right).
\end{equation}
In fact, following the lines of the proof of Lemma \ref{lema6}, it can be seen that
\begin{equation}\label{sesentaydos}
\mathbb{E}\frac{1}{S_{n-2}(x)+V}\leq 2\log \frac{27}{16}\mathbb{E}\frac{1}{S_{n-2}(x)+1}+\left(\log 4-2\log \frac{27}{16}\right)P(S_{n-2}(x)=0).
\end{equation}
On the other hand, it follows from (\ref{treintayuno}) that
\[
\mathbb{E}\frac{1}{S_{n-2}(x)+1}=\frac{1}{(n-1)x}P(S_{n-1}(x)\geq 1)=\frac{1}{x}\int_0^x(1-u)^{n-2}du
\]
\[
\leq \int_0^1e^{-(n-2)x\theta}d\theta \leq e^{2b_n/n}\int_0^1e^{-nx\theta}d\theta,
\]
where the last inequality follows from (\ref{cincuentaysiete}). This and (\ref{sesentaydos}) show claim (\ref{sesentayuno}).

Finally, since $0<x\leq 1/2$, we see from (\ref{treceestrella}) that $S_{n-2}(1-x)\geq S_{n-2}(1/2)$. We thus have from (\ref{sesentayuno})
\[
\mathbb{E}\frac{1}{S_{n-2}(1-x)+V}\leq \mathbb{E}\frac{1}{S_{n-2}(1/2)+V}
\leq e^{2b_n/n}\left(2\log \frac{27}{16}\int_0^1e^{-n\theta /2}d\theta +\left(\log 4-2\log \frac{27}{16}\right)e^{-n/2}\right)
\]
\[
\leq e^{2b_n/n}\left(\frac{4}{n}\log \frac{27}{16} +e^{-n/2}\right),
\]
which, in conjunction with (\ref{sesenta}) and (\ref{sesentayuno}), shows the result for $m=1$.

For $m>1$, the proof follows by induction on $m$, taking into account (\ref{cincuentaycinco}), (\ref{cincuentayseis}), and the fact that for any $0\leq \theta\leq 1$ and $j\in \mathbb{N}_0$, we have
\[
\mathbb{E}e^{-n\alpha_j(\theta)W_n(x)}=\mathbb{E}e^{-\alpha_j(\theta)(S_{n-2}(x)+V)}=\mathbb{E}e^{-\alpha_j(\theta)S_{n-2}(x)}\mathbb{E}e^{-\alpha_j(\theta)V}\leq e^{2b_n/n}e^{-n\alpha_{j+1}(\theta)x}\left(\frac{\alpha_{j+1}(\theta)}{\alpha_j(\theta)}\right)^2,
\]
as follows from (\ref{sesentaytres})-(\ref{sesentaycinco}). This completes the proof.
\hfill \break \rightline{$\square$}

Denote by
\begin{equation}\label{sesentayseis}
J_n(m,a)=\sup _{x\in \overline{R}_n(a)}\mathbb{E}\frac{\varphi^2(x)}{\varphi^2\left(W_n^{(m)}(x)\right)},\quad m\in \mathbb{N}.
\end{equation}

\begin{lemma}\label{lema18}
Let $m\in \mathbb{N}$. Assume that $b_n\leq 1/\alpha_{m-1}(1)$. Then,
\begin{align*}
&J_n(m,a)
\cr
& \leq b_ne^{2b_nm/n}\left(2\log \frac{27}{16}\int_0^1\left(\frac{\alpha_{m-1}(\theta)}{\theta}\right)^2e^{-b_n\alpha_{m-1}(\theta)}d\theta
+\left(\log 4-2\log \frac{27}{16}\right)\alpha_{m-1}^2(1)e^{-b_n\alpha_{m-1}(1)}+\epsilon _n\right),
\end{align*}
where $b_n$ and $\epsilon_n$ are defined in (\ref{cincuentaysiete}) and (\ref{cincuentaynueve}), respectively.
\end{lemma}
\textit{Proof}.\  Set $\lambda=nx$. The result readily follows from Lemma \ref{lema17} after observing the following. First, (\ref{cincuentaysiete}) implies that $\lambda <b_n$. Second, the function $h(\lambda)=\lambda e^{-\lambda \alpha_{m-1}(1)}$, $\lambda \geq 0$, increases in $\lambda < b_n\leq 1/\alpha_{m-1}(1)$.
\hfill \break \rightline{$\square$}

The following is the main result of this section.

\begin{theorem}\label{th19}
    Let $m\in \mathbb{N}$ and let $b_n$ be as in (\ref{cincuentaysiete}). Suppose that $b_n\leq 1/\alpha_{i-1}(1)$, for some $1\leq i \leq m$. Under assumption ($H_2$), we have
    \[
    \frac{\left\Vert \varphi^2g^{(2)}\right\Vert }{n}\left(1-J_n(m+1,a)\right)\leq \sqrt{2}\left(i+\sum_{k=i}^mJ_n(k,a)\right)\Vert B_nf-f\Vert ,
    \]
    where $J_n(m,a)$ is defined in (\ref{sesentayseis}).
\end{theorem}
\textit{Proof}.\  Our starting point is identity (\ref{cincuentaycuatro}). Let $k=0,1,\ldots ,i-1$. By the second equality in (\ref{catorce}), we have
\[
\left\vert \varphi^2(x)\left(\left(B_n^{k+1}g\right)^{(2)}(x)-\left(B_n^kg\right)^{(2)}(x) \right)\right\vert =
\left\vert \varphi^2(x)\left( B_n \left( B_n^{k+1}f-B_n^kf   \right)\right)^{(2)}(x) \right\vert
\]
\[
=\frac{2}{\varphi^2(x)}\left\vert \mathbb{E}\left(B_n^{k+1}f-B_n^kf\right)\left(\frac{S_n(x)}{n}\right)K_2(x;S_n(x))\right\vert
\]
\begin{equation}\label{sesentaysiete}
\leq \Vert B_nf-f\Vert  \frac{2}{\varphi^2(x)} \mathbb{E}\left\vert K_2(x;S_n(x))\right\vert \leq n\sqrt{2} \Vert B_nf-f\Vert ,
\end{equation}
where the last inequality follows from (\ref{once}) and Schwarz's inequality.

Let $k=i,\ldots ,m$ and set $h=B_ng-g$. By Lemma \ref{lema15} and (\ref{sesentayseis}), we get
\[
\left\vert \varphi^2(x)\left(\left(B_n^{k+1}g\right)^{(2)}(x)-\left(B_n^kg\right)^{(2)}(x) \right)\right\vert =
\left\vert \varphi^2(x) \left(B_n^kh\right)^{(2)}(x) \right\vert
=\left(\frac{n-1}{n}\right)^k\left\vert \varphi^2(x)\mathbb{E}h^{(2)}\left(W_n^{(k)}(x)\right)\right\vert
\]
\begin{equation}\label{sesentayocho}
\leq \left\Vert \varphi^2(B_ng-g)^{(2)}\right\Vert J_n(k,a)
\leq n\sqrt{2}\Vert B_nf-f\Vert J_n(k,a),
\end{equation}
where the last inequality follows from Lemma \ref{lema9}. Similarly,
\[
\left\vert \varphi^2(x)\left(B_n^{m+1}g\right)^{(2)}(x)\right\vert =\left(\frac{n-1}{n}\right)^{m+1}\left\vert \varphi^2(x)\mathbb{E}g^{(2)}\left(W_n^{(m+1)}(x)\right)\right\vert \leq \left\Vert \varphi^2g^{(2)}\right\Vert J_n(m+1,a).
\]
In view of (\ref{sesentaysiete}) and (\ref{sesentayocho}), the result follows from (\ref{cincuentaycuatro}) and assumption ($H_2$).
\hfill \break \rightline{$\square$}

Obviously, Theorem \ref{th19} makes sense only when $J_n(m+1,a)<1$. This can be achieved by choosing $m$ large enough, as follows from the upper bound in Lemma \ref{lema18} and the fact that $\alpha _m(\theta)\rightarrow 0$, as $m\rightarrow \infty$, $0<\theta \leq 1$.

\section{Proof of Theorem \ref{th1}}\label{sec6}

\underline{Upper bound}.
Let $a>0$. From Theorem \ref{th3}, we have
\begin{equation}\label{sesentaynueve}
\varlimsup _{n\rightarrow \infty}\sup_{f\in \mathfrak{F}}\frac{\omega_2^{\varphi}\left(f;1/\sqrt{n}\right)}{\|B_nf-f\|} \leq  4+\log 4 \varlimsup _{n\rightarrow \infty}\sup_{f\in \mathfrak{F}}\frac{n^{-1}\left\Vert \varphi^2(B_nf)^{(2)}\right\Vert }{\|B_nf-f\|}.
\end{equation}
Assume ($H_1$). Applying Theorems \ref{th4} and \ref{th14}, the right-hand side in (\ref{sesentaynueve}) is bounded above by
\begin{equation}\label{setenta}
4+\frac{\sqrt{2}\left(\sqrt{2}+1\right)}{1-0.99 K(a)/3}\log 4 ,
\end{equation}
where $K(a)$ is defined in (\ref{cincuenta}).

Assume ($H_2$). Let $m\in \mathbb{N}$ and $k=1,\ldots ,m$. Define
\begin{equation}\label{setentayuno}
J(k,a)=a\left(2 L_k(a)\log \frac{27}{16}+\left(\log 4-2\log \frac{27}{16}\right)\alpha_{k-1}^2(1)e^{-a\alpha_{k-1}(1)}\right),
\end{equation}
where
\begin{equation}\label{setentaydos}
L_k(a)=\int_0^1\left(\frac{\alpha_{k-1}(\theta)}{\theta}\right)^2e^{-a\alpha_{k-1}(\theta)}d\theta.
\end{equation}
By (\ref{cincuentaysiete}) and Lemma \ref{lema18}, we see that
\[
\varlimsup _{n\rightarrow \infty}J_n(k,a)\leq J(k,a).
\]
Therefore, by Theorem \ref{th19}, the right-hand side in (\ref{sesentaynueve}) can be bounded above by
\begin{equation}\label{setentaytres}
4+\frac{\sqrt{2}\left(i+\sum_{k=i}^mJ(k,a)\right)}{1-J(m+1,a)}\log 4,
\end{equation}
where $i\in \{1,\ldots ,m\}$ is the first integer such that $a<1/\alpha_{i-1}(1)$.

After numerically experimenting with expressions (\ref{setenta})-(\ref{setentaytres}), with the aid of Mathematica 10.2, we choose $a=7.2$, $m=20$, and $i=13$. Notice that this software system  allows to handle the involved iterates $\alpha _k(\theta)$, $k=0,\ldots , m$, defined in (\ref{sesentaytres}), and evaluate the integrals $L_k(a)$ with the command \lq NIntegrate\rq . As a result, both expressions (\ref{setenta}) and (\ref{setentaytres}) turn out to be less than $74.8$. This shows the upper bound in Theorem \ref{th1}.

\underline{Lower bound}.
For each $n\in \mathbb{N}$, define the function $f_n\in C[0,1]$ as follows
\begin{equation}\label{AA}
f_n\left(\frac{2-\sqrt{2}}{n}\right)=1,\quad f_n\left(\frac{1}{n}\right)=-0.8,\quad f_n\left(\frac{2}{n}\right)=-1,\quad f_n\left(\frac{3}{n}\right)=0.04,\quad f_n\left(\frac{2+\sqrt{2} }{n}\right)=1.
\end{equation}
If $y\in [(2-\sqrt{2})/n,(2+\sqrt{2})/n]$, then $f_n(y)$ is piecewise linear interpolating the values in (\ref{AA}), whereas if $y\in [0,1]\setminus [(2-\sqrt{2})/n,(2+\sqrt{2})/n]$, then $f_n$ is constantly $1$.

Taking the weighted second order difference of $f$ at the point $2/n$, it can be checked that
\begin{equation}\label{BB}
\lim _{n\rightarrow \infty}\omega _2^{\varphi}\left(f_n;\frac{1}{\sqrt{n}}\right)=4.
\end{equation}

On the other hand, set $\lambda=nx$ and denote $g(\lambda)=f_n(\lambda /n)$, $0\leq \lambda \leq n$. It turns out that
\begin{equation}\label{CC}
B_nf_n(x)-f_n(x)=\mathbb{E}f_n\left(\frac{S_n(x)}{n}\right)-f_n\left(\frac{\lambda}{n}\right)=\mathbb{E}g\left(S_n\left(\frac{\lambda}{n}\right)\right)-g(\lambda),\quad 0\leq \lambda \leq n.
\end{equation}

Let $\lambda _0>2+\sqrt{2}$ and assume that $n\geq \lambda _0$ is large enough. Using (\ref{treceestrella}), (\ref{treintaycuatro}), and (\ref{treintaysiete}), we have for $\lambda > \lambda _0$
\[
\left\vert \mathbb{E}g\left(S_n\left(\lambda /n\right)\right)-g(\lambda)\right\vert
=\left\vert \mathbb{E}\left(g\left(S_n\left(\lambda /n\right)\right)-g(\lambda)\right)1_{\left\{S_n(\lambda/n)\leq 2+\sqrt{2}\right\}} \right\vert \leq 2P\left(S_n\left(\lambda /n\right)\leq 3\right)
\]
\begin{equation}\label{DD}
 \leq 2P\left(S_n\left(\lambda _0 /n\right)\leq 3\right)   \leq  2P\left(N_{\lambda_0}\leq 3\right)+2d_{TV}\left(S_n\left(\lambda_0 /n\right),N_{\lambda_0}\right)\leq  2P\left(N_{\lambda_0}\leq 3\right)+2\frac{D(\lambda_0)}{n}.
\end{equation}

Recalling (\ref{AA}), denote by
\begin{align*}
G(\lambda)=\mathbb{E}g(N_{\lambda})&=P(N_{\lambda}=0)-0.8P(N_{\lambda}=1) \\
& -P(N_{\lambda}=2)+0.04P(N_{\lambda}=3)+1-P(N_{\lambda}\leq 3),\qquad 0\leq \lambda \leq \lambda_0.
\end{align*}
By (\ref{cuarentaydosestrella}), (\ref{treintaycuatro}), and (\ref{treintaysiete}), we have for $\lambda \leq \lambda _0$
\[
\left\vert \mathbb{E}g\left(S_n\left(\lambda /n\right)\right)-g(\lambda)\right\vert
\leq \left\vert \mathbb{E}g(N_{\lambda})-g(\lambda)\right\vert +\left\vert \mathbb{E}g\left(S_n\left(\lambda /n\right)\right)-\mathbb{E}g(N_{\lambda})\right\vert
\]
\begin{equation}\label{EE}
\leq \left\vert G(\lambda)-g(\lambda) \right\vert
+d_{TV}(S_n(\lambda /n),N_{\lambda })  \leq \left\vert G(\lambda)-g(\lambda) \right\vert  +\frac{D(\lambda_0)}{n}.
\end{equation}

Choose $\lambda_0$ and $n$ large enough so that the upper bound in (\ref{DD}) is less than $1/2$, say. On the other hand, numerical computations, performed again with Mathematica 10.2, show that
\[
\sup _{0\leq \lambda \leq \lambda _0}\left\vert G(\lambda)-g(\lambda)\right\vert \leq 0.79 \ldots
\]
We therefore have from (\ref{CC}), (\ref{DD}), and (\ref{EE})
\[
\left\Vert B_nf_n-f_n\right\Vert \leq 0.8,
\]
for large enough $n$. By (\ref{BB}), this implies that
\[
\varlimsup_{n\rightarrow \infty }\sup _{f\in \mathfrak{F}}\frac{\omega_2^{\varphi}\left(f;1/\sqrt{n}\right)}{\|B_nf-f\|} \geq \varlimsup_{n\rightarrow \infty }\ \frac{\omega_2^{\varphi}\left(f_n;1/\sqrt{n}\right)}{\|B_nf_n-f_n\|} \geq \frac{4}{0.8}=5.
\]
This shows the lower bound in Theorem \ref{th1}, and concludes the proof.
\hfill \break \rightline{$\square$}

\section{Proof of Theorem \ref{th2}}\label{sec7}

Let $f\in \mathfrak{G}$. By the first equality in (\ref{catorce}), we have
\begin{equation}\label{CCCCC}
\left\vert \varphi^2(x)(B_nf)^{(2)}(x)\right\vert \leq n^2\varphi^2(x)\omega_2 \left( f;1/n\right).
\end{equation}
Let $a>0$. Define
\[
N_1:=N_1(a,f)=\left\{n\in \mathbb{N}:\left\Vert \varphi^2 (B_nf)^{(2)}\right\Vert \text{ is attained on }\overline{R}_n(a)\right\},\qquad N_2:=\mathbb{N}\setminus N_1
\]
Suppose that $N_1$ is not finite. By (\ref{AAAAA}) and (\ref{CCCCC}), we see that
\[
\lim_{\substack{n\to \infty \\ n\in N_1}}\left\Vert \varphi^2 (B_nf)^{(2)}\right\Vert \leq a \lim_{\substack{n\to \infty \\ n\in N_1}}    n\ \omega_2 \left(f;1/n\right)=0.
\]
We thus have from Theorem \ref{th3}
\begin{equation}\label{DDDDD}
\varlimsup_{\substack{n\to \infty \\ n\in N_1}} \frac{\omega_2^{\varphi}\left(f;1/\sqrt{n}\right)}{\|B_nf-f\|}\leq 4+ \log 4\varlimsup_{\substack{n\to \infty \\ n\in N_1}}\  \frac{\left\Vert \varphi ^2(B_nf)^{(2)}\right\Vert}{n\|B_nf-f\|}=4,
\end{equation}
where we have used that $f\in \mathfrak{F}$ together with the well known fact that $\Vert B_nf-f\Vert =o(1/n)$ if and only if $f$ is affine.

Suppose that $N_2$ is not finite. By (\ref{setenta}), we have
\begin{equation}\label{DDD}
\varlimsup_{\substack{n\to \infty \\ n\in N_2}} \frac{\omega_2^{\varphi}\left(f;1/\sqrt{n}\right)}{\|B_nf-f\|}\leq 4+\frac{\sqrt{2}\left(\sqrt{2}+1\right)}{1-0.99 K(a)/3}\log 4 .
\end{equation}
Since $a$ is arbitrary, and $k(s)\rightarrow \sqrt{3}$, as $s\to \infty$ (see (\ref{cincuenta})), inequality (\ref{BBBBB}) follows from (\ref{DDDDD}) and (\ref{DDD}).

Finally, let $f\in \mathfrak{F}\cap C^1[0,1]$ and $0 <h \leq \delta $. We have from (\ref{dos})
\[
\vert f(x+h)-2f(x)+f(x-h)\vert = h \left\vert \mathbb{E}\left(f^{(1)}(x+h\beta _1)-f^{(1)}(x-h\beta _1)\right)   \right\vert \leq h\omega _1\left(f^{(1)};2h\right).
\]
By (\ref{AAAAA}), this readily implies that $f\in \mathfrak{G}$. This concludes the proof.
\hfill \break \rightline{$\square$}

\textbf{Funding}: The first author is supported by Research Project DGA (E48\_23R). The second author is supported by Junta de Andaluc\'\i a (Research Group FQM-0178).

\bibliographystyle{unsrtnat}

\begin{thebibliography}{99}


\bibitem{abelleviatanrasa}
Abel, U., Leviatan, D., Ra\c{s}a, I.: On the $q$-monotonicity preservation of Durrmeyer-type operators. Mediterr. J. Math.  (2021), \textbf{18}:173.

\bibitem{acurasasteopoaie}
Acu, A.M., Ra\c{s}a, I., \c{S}teopoaie, A.E.: Bernstein–Kantorovich operators, approximation and shape preserving properties. Rev. Real Acad. Cienc. Exactas Fis. Nat. Ser. A-Mat. (2024), \textbf{118}:107.

\bibitem{adellanoz}
Adell, J.A., Anoz, J.M.: Signed binomial approximation of binomial mixtures via differential calculus for linear operators. J. Statist. Plann. Inference \textbf{138} (2008), 3687--3695.

\bibitem{adellanozlekuona}
Adell, J.A., Anoz, J.M., Lekuona, A.: Exact values and sharp estimates for the total variation distance between binomial and Poisson distributions. Adv. Appl. Prob. \textbf{40} (2008), 1033--1047.

\bibitem{rm2019}
Adell, J.A., C\'{a}rdenas-Morales, D.: On the $10$th central moment of the Bernstein polynomials. Results Math. (2019),  \textbf{74}:113.

\bibitem{rm2022}
Adell, J.A., C\'{a}rdenas-Morales, D.: Asymptotic and non-asymptotic results in the approximation by Bernstein polynomials. Results Math. (2022), \textbf{77}:166.

\bibitem{adelllekuona}
Adell, J.A., Lekuona, A.: Binomial convolution and transformations of Appell polynomials. J. Math. Anal. Appl. \textbf{456} (2005), 16--33.

\bibitem{adellsanguesa}
Adell, J.A., Sang\"{u}esa, C.: A strong converse inequality for gamma-type operators. Constr. Approx. \textbf{15} (1999), 537--551.

\bibitem{barbourhall}
Barbour, A.D., Hall, P.: On the rate of Poisson convergence. Math. Proc. Cambridge Philos. Soc. \textbf{95}(3) (1984), 473--480.

\bibitem{bustamante}
Bustamante, J.: Estimates of positive linear operators in terms of second order moduli. J. Math. Anal. Appl. \textbf{345} (2008), 203--212.

\bibitem{bustamantee}
Bustamante, J.: Baskakov-Kantorovich operators reproducing affine functions: inverse results. J. Numer. Anal. Approx. Theory \textbf{51}(1) (2022), 67--82.

\bibitem{bustamantequesada}
Bustamante, J., Quesada, J.M.: A property of Ditzian-Totik second order moduli. Appl. Math. Lett. \textbf{23} (2010), 576--580.

\bibitem{chenditzian}
Chen, W., Ditzian, Z.: Strong converse inequality for Kantorovich polynomials. Constr. Approx. \textbf{10} (1994), 107--129.

\bibitem{chihara}
Chihara, T.S.: An Introduction to Orthogonal Polynomials. Gordon and Breach, New York (1978).

\bibitem{deheuvelspfeifer}
Deheuvels, P., Pfeifer, D., Puri, M.L.: A new semigroup technique in Poisson approximation. Semigroup Forum \textbf{38}(2) (1989), 189--201.

\bibitem{dellavecchia}
Della Vecchia, B.: Direct and converse results by rational operators. Constr. Approx. \textbf{12} (1996), 271--285.

\bibitem{ditzianivanov}
Ditzian, Z., Ivanov, K.G.: Strong converse inequalities. J. Anal. Math. \textbf{61} (1993), 61--111.

\bibitem{ditziantotik}
Ditzian, Z., Totik, V.: Moduli of Smoothness. Springer-Verlag, New York (1987).

\bibitem{finta}
Finta, Z.: Direct and converse results for $q$-Bernstein operators. Proc. Edinb. Math. Soc. \textbf{52}(2) (2009), 339--349.

\bibitem{gadjev}
Gadjev, I.: Strong converse result for uniform approximation by Meyer-K\"{o}nig and Zeller operator. J. Math. Anal. Appl. \textbf{428} (2015), 32--42.

\bibitem{gavreaetall}
Gavrea, I., Gonska, H.H., P\u{a}lt\u{a}nea, R., Tachev, G.: General estimates for the Ditzian-Totik modulus. East J. Approx. \textbf{9}(2) (2003), 175--194.

\bibitem{guoqi}
Guo, SH., Qi, Q.: Strong converse inequalities for Baskakov operators. J. Aprox. Theory \textbf{124}(2) (2003), 219--231.

\bibitem{knoopzhou1994}
Knoop, H.B., Zhou, X.L.: The lower estimate for linear positive operators, II. Results Math. \textbf{25} (1994),  315--330.

\bibitem{knoopzhou1995}
Knoop, H.B., Zhou, X.L.: The lower estimate for linear positive operators, I. Constr. Approx. \textbf{11} (1995), 53--66.

\bibitem{lopezsalamanca}
L\'{o}pez-Bl\'{a}zquez, F., Salamanca Mi\~{n}o, B.: Binomial approximation to hypergeometric probabilities. J. Statist. Plann. Inference \textbf{87} (2000), 21--29.

\bibitem{paltanea2004}
P\u{a}lt\u{a}nea, R.: Approximation Theory Using Positive Linear Operators. Birkh\"auser Boston, Inc., Boston, MA (2004).

\bibitem{paltanea2018}
P\u{a}lt\u{a}nea, R.: Asymptotic constant in approximation of twice differentiable functions by a class of positive linear operators. Results Math. \textbf{73}(2) (2018), paper no. 64, 10 pp.

\bibitem{roos}
Roos, B.: Binomial approximation to the Poisson binomial distribution. The Krawtchouk expansion. Theory Probab. Appl. \textbf{45} (2000), 258--272.

\bibitem{sanguesa}
Sang\"{u}esa, C.: Lower estimates for centered Bernstein-type operators. Constr. Approx. \textbf{18}(1) (2002), 145--159.

\bibitem{totik}
Totik, V.: Strong converse inequalities. J. Approx. Theory \textbf{76} (1994), 369--375.

\bibitem{totikk}
Totik, V.: Approximation by Bernstein polynomials. Amer. J. Math. \textbf{116} (1994), 995--1018.


\end{thebibliography}

\end{document}